%% file: strong912.tex
\documentclass[11pt]{amsart} 
\usepackage{amsfonts,amscd,amssymb}
\usepackage[dvips]{epsfig,graphics}
\usepackage{graphicx}

\newtheorem{thm}{Theorem}[section]
\newtheorem{lemma}[thm]{Lemma}
 \newtheorem{prop}[thm]{Proposition}

\newtheorem{defn}[thm]{Definition}

\newcommand{\R}{\mathbb R}

\renewcommand{\o}{\operatorname}

\newcommand{\oto}{\o{O}(2,1)} \newcommand{\soto}{{\o{SO}(2,1)^0}}

\newcommand{\rto}{\R^{2,1}}

\newcommand{\E}{{\mathbb A}^{2,1}}
\newcommand{\affine}{\E}
\newcommand{\isoaff}{\o{Aff}(\E)}

\newcommand{\hyp}{{\mathbb H}^2}
\newcommand{\psl}{\o{Isom}(\hyp)}
\newcommand{\Hyp}{{\mathcal H}}

\newcommand{\vx}{{\mathsf x}}
\newcommand{\vy}{{\mathsf y}}

\newcommand{\vu}{{\mathsf u}}
\newcommand{\vv}{{\mathsf v}}
\newcommand{\vw}{{\mathsf w}}

\newcommand{\xo}[1]{{\mathsf x}^0(#1)}
\newcommand{\xp}[1]{{\mathsf x}^+(#1)}
\newcommand{\xm}[1]{{\mathsf x}^-(#1)}
\newcommand{\xpm}[1]{{\mathsf x}^{\pm}(#1)}

\newcommand{\ldot}[1]{\mathbb B(#1)}

\newcommand{\cqfd}{\begin{flushright}$\square$\end{flushright}}

\newcommand{\vg}{\vv_g}
\newcommand{\vh}{\vv_h}

\newcommand{\lag}{\lambda_g}
\newcommand{\lah}{\lambda_h}

\newcommand{\lagp}{\lambda_{\gp}}
\newcommand{\lahp}{\lambda_{\hp}}

\newcommand{\ga}{\alpha}
\renewcommand{\gg}{\gamma}
\newcommand{\gG}{\Gamma}
\newcommand{\gh}{\eta}
\newcommand{\gw}{\omega}
\newcommand{\gl}{\lambda}
\newcommand{\gk}{\kappa}
\newcommand{\cg}{C_\gg}

\newcommand{\ggp}{\gg_\phi}
\newcommand{\ghp}{\gh_\phi}
\newcommand{\gkp}{\gk_\phi}
\newcommand{\gp}{g_\phi}
\newcommand{\hp}{h_\phi}

\renewcommand{\R}{\mathbb R}

\newcommand{\B}{\mathbb B}
\newcommand{\LL}{\mathbb L}

\newcommand{\G}{\mathfrak G}
\newcommand{\g}{\mathfrak g}
\newcommand{\U}{\mathcal U}

\begin{document}

\title{Strong marked isospectrality of affine Lorentzian groups}

\author{Virginie Charette}
     \address{Department of Mathematics\\ University of Manitoba\\
     Winnipeg, Manitoba, Canada}
     \email{charette@cc.umanitoba.ca}
\author{Todd Drumm}
     \address{Department of Mathematics\\ University of Pennsylvania\\
     Philadelphia, PA}
     \email{tad@math.upenn.edu}

\thanks{The authors would like to thank Bill Goldman for many helpful conversations.}

\date{\today}

\begin{abstract}
The Margulis invariant $\ga$ is a function on $H^1(\gG,\rto)$, where
$\gG$ is a group of Lorentzian transformations acting on
$\R^{2,1}$, that contains no elliptic elements. The spectrum of $\gG$ is the image of all  $\gg\in\gG\setminus Id$
under the map $\ga$.  If the underlying linear group of $\gG$ is
fixed, Drumm and Goldman proved that the spectrum defines the
translational part completely. In this note, we strengthen this
result by showing that isospectrality holds for any free product of cyclic 
groups of given rank, up to conjugation in the group of affine transformations of 
$\rto$, as long as it is non-radiant and that its linear part is discrete and 
non-elementary.  In particular, isospectrality holds when the linear
part is a Schottky group.
\end{abstract}

\maketitle

\section{Introduction}

The {\em Margulis invariant}, $\ga(\gg)$, of a fixed point
free affine hyperbolic transformation $\gg$ measures oriented translation along the unique $\gg$-invariant line in $\E$, affine 2+1
Minkowski space. This is an important invariant, which Margulis used in \cite{Ma1,Ma2} to prove the existence of proper
actions by purely hyperbolic affine groups on $\R^3$.  In a previous paper, we extended its definition to affine transformations with parabolic linear part~\cite{CD}.

Let $\Gamma_0\subset\soto$.  Since $\soto\cong\psl$, the group of isometries of the 
hyperbolic plane, an element is hyperbolic, parabolic or elliptic, depending on its 
fixed point set.  An {\em affine deformation} of
$\Gamma_0$ is an isomorphism into the group of affine transformations of $\E$: 
\begin{equation*}
     \phi:\Gamma_0\rightarrow\isoaff
\end{equation*}
such that $\LL\circ\phi=\Gamma_0$, where $\LL$ denotes projection
onto the linear part.  A {\em Schottky subgroup of $\soto$} is a
discrete, freely generated subgroup of $\soto$ whose non-identity
elements are hyperbolic.  Drumm-Goldman~\cite{DG} proved that, given
a fixed linear part $\Gamma_0\subset\soto$ that is Schottky, the marked Margulis spectrum
determines the affine deformation up to translational conjugacy:

\begin{thm}\label{thm:DG}{\em [Weak Isospectrality]}~(Drumm-Goldman)
Let $\Gamma_0$ be a discrete, purely hyperbolic, free subgroup of
$\soto$ and $\phi_1$, $\phi_2$, a pair of affine deformations of
$\Gamma_0$ such that, for every $\gg\in\Gamma_0$:
\begin{equation*}
\ga(\phi_1(\gg))=\ga(\phi_2(\gg)).
\end{equation*}
Then $\phi_1=T\circ\phi_2\circ T^{-1}$, where $T$ is a
translation.
\end{thm}
Kim~\cite{K} generalized this result to
higher dimensions.

In this paper, we generalize Theorem~\ref{thm:DG} by removing the
assumption that the linear part $\Gamma_0$ is known. We thus
consider isomorphisms from an abstract group $\G$, which is assumed to
be a finitely generated free product of cyclic groups, to $\isoaff$:
\begin{equation*}
\phi : \G\longrightarrow\isoaff .
\end{equation*}
A group $\Gamma\in\isoaff$ is called {\em radiant} if there is a
point $x\in\E$ fixed by $\Gamma$. By Drumm-Goldman's theorem, a group
freely generated by hyperbolic isometries is radiant if and only
if the Margulis invariant is identically zero on $\Gamma$.  We
will show that if $\phi(\G)$ is not radiant and if $\LL\circ\phi$ is an isomorphism onto a discrete, non-elementary subgroup, 
then the growth
of the Margulis invariant as a function of word length in the
group  determines the $\isoaff$-conjugacy class of $\phi(\G)$.  (Recall that a group $G\subset\psl$ is {\em elementary} if it admits a finite orbit.)

\begin{thm}\label{thm:Main}{\em [Strong isospectrality]}~
Let $\G$ be a finitely generated, free product of cyclic groups.  Let $\phi_1$ and $\phi_2$
be isomorphisms:
\begin{equation*}
     \phi_i : \G\longrightarrow\isoaff ,
\end{equation*}
such that $\phi_i(\G)$ is not radiant, and $\LL\circ\phi_i$ is an
isomorphism onto a discrete and non-elementary subgroup of $\soto$. Suppose that for every
$\g\in\G$ such that both $\phi_1(\g)$ and $\phi_2(\g)$ are hyperbolic:
\begin{equation*}
\ga(\phi_1(\g))=\ga(\phi_2(\g)).
\end{equation*}
Then $\phi_1(\G)$ and $\phi_2(\G)$ are $\isoaff$-conjugate.
\end{thm}
This includes the case where the linear parts of $\phi_1(\G)$ and
$\phi_2(\G)$ are both Schottky subgroups of $\soto$.

We do not require a priori that $\phi_1(\g)$ and $\phi_2(\g)$ be of
the same type 
(i.e. hyperbolic, parabolic or elliptic), given $\g\in\G$.   But our
assumptions do impose the following restrictions.  Since $\LL\circ\phi_i$ is an
isomorphism onto a discrete group, 
$\g\in\G$ is of finite order if and only if $\phi_i(\g)$ has elliptic
linear part.  
Furthermore, no two elements share common fixed points unless they are
powers of 
each other. Finally, the fact that $\phi_i(\G)$ is non-elementary implies the existence of a hyperbolic element whose fixed point set is not invariant under the action of any elliptic element in the generating set.

The outline of the proof is as follows.  Choose an appropriate
generating set for $\G$, so that the generating sets of both $\phi_1(\G)$ and $\phi_2(\G)$ consist entirely of hyperbolic elements.

Next, reduce to the case of rank two groups whose linear parts are
Schottky.  Choose appropriate
conjugates of $\phi_1(\G)$ and $\phi_2(\G)$ so that the eigensystem of one 
corresponding pair
of generators is the same, and the other corresponding pair of generators 
shares one eigendirection.

Equality of the Margulis invariant on hyperbolic words of the form $\gh^n\gg^m$ 
determines the remaining
eigendirections.  We then show that the eigenvalues are also the same, thus 
concluding that the
linear parts are equal.  

By Theorem~\ref{thm:DG}, the affine groups are 
thus equal up to
translational conjugacy.  We will include a proof of Theorem~\ref{thm:DG} here for the sake of completeness.

\section{Preliminaries}
Let $\affine$ denote three-dimensional affine space with the following 
additional structure.
Its associated vector space of directions
\begin{equation*}
     \R^{2,1}=\{ p-q\mid p,q\in\affine\} ,
\end{equation*}
which is isomorphic to $\R^3$ as a vector space, is endowed with the 
standard {\em Lorentzian
scalar product}:
\begin{equation*}
     \ldot{\vx,\vy}=x_1y_1+x_2y_2-x_3y_3,
\end{equation*}
where $\vx=[x_1~x_2~x_3]^T$, $\vy=[y_1~y_2~y_3]^T\in\R^3$. Thus, $\affine$ is
{\em Minkowski (2+1)-spacetime}.

A non-zero
vector $\vx$ is said to be {\em null} (resp. {\em timelike}, {\em
spacelike}) if $\ldot{\vx,\vx}=0$ (resp. $\ldot{\vx,\vx}<0$, 
$\ldot{\vx,\vx}>0$).
A null vector is {\em future-pointing} if its third coordinate is positive 
-- this
corresponds to choosing a connected component of the set of timelike 
vectors, or a
{\em time-orientation}.

The {\em Lorentzian cross-product} is the unique bilinear map:
\begin{equation*}
     \boxtimes :\R^{2,1}\times\R^{2,1}\longrightarrow\R^{2,1}
\end{equation*}
such that $\ldot{\vu,\vv\boxtimes\vw}=\mbox{Det}(\left[\vu~
\vv~ \vw\right]) $.  It has the following properties:
\begin{itemize}
\item $\ldot{\vu,\vv\boxtimes\vw}=\ldot{\vv,\vw\boxtimes\vu}$;
\item $\ldot{\vv,\vv\boxtimes\vw}=0$;
\item 
$\ldot{\vv\boxtimes\vw,\vv\boxtimes\vw}=\ldot{\vv,\vw}^2-\ldot{\vv,\vv}\ldot{\vw,\vw}$.
\end{itemize}

\subsection{Affine deformations of a linear group}

Let $\isoaff$ denote the group of all affine transformations that
preserve the Lorentzian scalar product on the space of directions.
Choosing an origin in $\affine$ allows one to write an affine
transformation as the composition of a linear transformation with a translation: 

\begin{equation*}
\gg (x)  =  g(x)+\vg;
\end{equation*}
$g$ is the {\em linear part} of $\gg$ and $\vg$ its {\em translational part}.

Thus $\isoaff$ is isomorphic to $\oto\ltimes\R^{2,1}$.
Denote projection onto the linear part of an affine transformation by:
\begin{equation*}
 \LL:\isoaff\rightarrow\oto .
\end{equation*}

Let $\G$ be a finitely generated group and $\phi:\G\longrightarrow\isoaff$ an isomorphism; denote the linear part of $\phi$ by $\Phi=\LL\circ\phi$.  We call $\phi$ an {\em affine deformation} of $\Phi$.

Fix the linear part $G_\Phi=\Phi(\G)$.  The left action induces a $G_\Phi$-module structure on $\rto$, which we will denote by $V_\Phi$.  Now $\phi$ is a homomorphism if and only if the translational part satisfies the cocycle condition:

\begin{equation*}
\vv_{gh}=\vg+g(\vh),
\end{equation*}
for every $g,h\in G_\Phi$.  Thus an affine deformation of $\Phi$ may be interpreted as a cocycle $u\in Z^1(G_\Phi,V_\Phi)$, where $u(g)=\vg$.  So we write:

\begin{equation*}
\phi =(\Phi,u).
\end{equation*}

Suppose now that $\psi =\tau\phi\tau^{-1}$, where $\tau\in\isoaff$ is the pure translation: $\tau(x)=x+\vv$, $\vv\in\rto$.  Then the linear part of $\psi$ is $\Phi$ and for $\g\in\G$, the translational part of $\psi(\g)$ is $u(\phi(\g))+\vv-\Phi(\g)(\vv)$.  In other words, the translational parts of $\phi(\g)$ and $\psi(\g)$ differ by the coboundary $\vv-\Phi(\g)(\vv)\in B^1(G_\Phi,V_\Phi)$.  Thus $H^1(G_\Phi,V_\Phi)$ corresponds to translational conjugacy classes of affine deformations of $\Phi$.

\subsection{Hyperbolic transformations}

We shall restrict our attention to those transformations whose
linear parts are in $\soto$, thus preserving orientation and
time-orientation. The isomorphism between $\soto$ and $\psl$ gives
rise to the following terminology.

 \begin{defn}
 Let $g\in\soto$ be
a non-identity element; 
\begin{itemize}
\item $g$ is {\em hyperbolic} if it has three,
distinct, real eigenvalues;
\item $g$ is {\em parabolic} if its only real eigenvalue is 1, corresponding to a one-dimensional null eigenspace;
\item $g$ is {\em elliptic} if its only real eigenvalue is 1, corresponding to a one-dimensional timelike eigenspace.
\end{itemize}

We also call $\gamma\in\isoaff$ {\em hyperbolic} (resp. {\em parabolic}, {\em elliptic}), if its linear
part $\LL(\gg)$ is hyperbolic (resp. parabolic, elliptic).   Denote by $\Hyp$ the set of
hyperbolic affine transformations.
\end{defn}

Suppose $g\in\soto$ is hyperbolic, with eigenvalues $\lag ,1,1/\lag$,
for some $0<\lag<1$.  Associated to $g$ is a {\em null frame}:
\begin{equation*}
\{ \xo{g},\xm{g},\xp{g}\} ,
\end{equation*}
where $\xm{g}$ (resp. $\xp{g}$) is a null, future-pointing $\lag$-eigenvector
(resp. $1/\lag$-eigenvector), chosen to be of unit Euclidean length, and
$\xo{g}$ is the unique 1-eigenvector such that $\ldot{\xo{g},\xo{g}}=1$ and
$\{ \xo{g},\xm{g},\xp{g}\}$ is positively oriented.  If $\gamma\in\Hyp$, we
write:
\begin{equation*}
\{ \xo{\gamma},\xm{\gamma},\xp{\gamma}\} := \{ 
\xo{\LL(\gg)},\xm{\LL(\gg)},\xp{\LL(\gg)}\}.
\end{equation*}

Observe that $\ldot{\xo{g},\xpm{g}}=0$ and thus $\xo{g}$ is a
positive scalar multiple of $\xm{g}\boxtimes\xp{g}$.  In fact,
since
$\ldot{\xm{g}\boxtimes\xp{g},\xm{g}\boxtimes\xp{g}}=\ldot{\xm{g},\xp{g}}^2$:
\begin{equation}\label{eq:xog}
\xo{g}=\frac{-1}{\ldot{\xm{g},\xp{g}}}\xm{g}\boxtimes\xp{g}.
\end{equation}

The following fact will also prove useful.
\begin{lemma}\label{lem:xpxo}
Let $g\in\soto$ be hyperbolic.   Then
$\xo{g}\boxtimes\xp{g}=\xp{g}$ and $\xm{g}\boxtimes\xo{g}=\xm{g}$.
\end{lemma}
\begin{proof}
Since the vector $\xp{g}$ is Lorentz--perpendicular to both itself and
$\xo{g}$, we know that $\xo{g}\boxtimes\xp{g} = k\xp{g}$. Taking
the Lorentzian inner product of both sides of the equation with
$\xm{g}$ we get: 
\begin{equation*}
\ldot{\xo{g}\boxtimes\xp{g},\xm{g}} =k\ldot{\xp{g},\xm{g}}.
\end{equation*}
The left hand side can be rearranged via a property of $\boxtimes$
described above: 
\begin{equation*}
\ldot{\xp{g}\boxtimes\xm{g},\xo{g}}=k\ldot{\xp{g},\xm{g}}.
\end{equation*}
Rewriting $\xp{g}\boxtimes\xm{g}$, using \eqref{eq:xog}, we obtain the
following:
\begin{equation*}
\ldot{\xm{g},\xp{g}}\ldot{\xo{g},\xo{g}}=k\ldot{\xp{g},\xm{g}}.
\end{equation*}
This implies that $k=1$.  The proof for $\xm{g}$ is similar.
\end{proof}

\subsection{The Margulis invariant.}
Every affine hyperbolic $\gg$ admits a unique spacelike line that is $\gamma$--invariant, denoted $C_{\gg}$.  Furthermore, $\cg$ is parallel to $\xo{\gg}$ and $\gg$ acts by translation
on $C_{\gg}$. On the subset $\Hyp\subset\isoaff$ of hyperbolic elements,
we define the {\em Margulis invariant} of $\gamma$ to be the
function:  
\begin{equation*}
\ga : \Hyp  \rightarrow \R 
\end{equation*}
such that: 
\begin{equation} \label{eq:alphadef}
\ga(\gg) = \B(\gg(x)-x,\xo{\gg}),
\end{equation}
where $x$ is an arbitrary point on $C_{\gg}$, the unique $\ga$-invariant line.

The action on $C_{\gg}$ is given by
\begin{equation*}
\gg(x)=x+\ga(\gg)\xo{\gg},
\end{equation*}
for every $x\in\cg$.

The following are elementary consequences of the definition.

\begin{lemma}\label{lem:hyperbolicfacts}{\em [Properties of $\ga$]}~
Suppose $\gg\in\Hyp$;
\begin{enumerate}
\item for any $x\in\affine$, 
$\ga(\gg)=\ldot{\gg(x)-x,\xo{\gg}}$;\label{hfact:pointinvariance}
\item $\ga(\gg)\neq 0$ if and only if $\gg$ acts 
freely;\label{hfact:fixedpoint}
\item for any $\eta\in\isoaff$, 
$\ga(\eta\gg\eta^{-1})=\ga(\gg)$;\label{hfact:conjugation}
\item for any $\gg\in\isoaff$, $\ga(\gg^{n})= |n|\ga(\gg)$. \label{hfact:power}
\end{enumerate}
\end{lemma}
\cqfd

Item~(\ref{hfact:fixedpoint}) implies that if $\phi=(\Phi,u)$ is radiant then the Margulis invariant is identically zero on the group (Theorem~\ref{thm:DG} proves the converse).  Any two hyperbolic elements in $\soto$ with the same eigenvalues, or trace, 
are conjugate.
Furthermore, it can be shown that any two hyperbolic elements in $\isoaff$ 
with the same
linear part and Margulis invariant are conjugate by a translation. Thus, 
the trace of
the linear part and the value of $\ga$ determine the conjugacy classes of 
hyperbolic
elements in $\isoaff$.

Since $\xo{\gg^{-1}}=-\xo{\gg}$, we have that 
$\alpha(\gg^{-1})=\ga(\gg)$.  When $\gg$
acts freely, the sign of $\ga(\gg)$ indicates the direction in which $\gg$ 
displaces
points on $\cg$.

Suppose that $\G$ is a finitely generated group and let $\Phi:\G\longrightarrow \soto$ be an isomorphism.  Keep in mind that we will be reducing to the case where $\Phi(\g)$ is hyperbolic for every $\g\in\G\setminus Id$.   So we can safely ignore non-hyperbolic elements in the group.  Fix an ordering on the hyperbolic elements of $\G\setminus Id$, $\left[\g_i\right]_{i\in {\mathbb N}}$.

The {\em marked Margulis spectrum} is determined as follows.  By Item~(\ref{hfact:conjugation}) of Lemma~\ref{lem:hyperbolicfacts}, two affine deformations of $\Phi$ yield the same values of $\alpha$ if they are translationally conjugate.   Define the following function:

\begin{align*}
 \ga_\Phi:H^1(G_\Phi,V_\Phi) &\longrightarrow \R^\G \\
[u] & \longmapsto \left[\ldot{u(\Phi(\g_i)),\xo{\Phi(\g_i)}}\right]_{i\in {\mathbb N}} .
\end{align*}
 
The marked Margulis spectrum of an affine deformation $\phi=(\Phi,u)$ is the image $\ga_\Phi([u])$.

\section{Step one: Reduction to a simple case}

In this section, we show how to reduce Theorem~\ref{thm:Main} to a simpler 
case.  Indeed, we may assume that our groups sharing a marked Margulis 
spectrum are rank-two groups generated by hyperbolic isometries and that their generators admit certain common 
characteristics.

Let $\{\g_1,\ldots\g_n\}$ be a generating set for $\G$, with no
relations except, possibly, $\g_j^{m_j}=Id$, for some of the $j$'s.   
Our first concern is for the linear parts of $\phi_1$ and $\phi_2$, so
it will be 
simpler for now to think of $\phi_1$ and $\phi_2$ as representations
into $\psl$.  
Thus a hyperbolic element has two fixed points on the boundary of the
hyperbolic 
plane, a parabolic element has one fixed point on the boundary of the
hyperbolic 
plane, and an elliptic element fixes a point inside the hyperbolic plane.  Two elements $g,h\in\psl$ commute if and only if they share the same set of fixed points.

Since $\phi_1(\G)$ is non-elementary, it contains a hyperbolic element $\g$; we may choose it such that its null eigenspace is not invariant  under the action of any of the elliptic elements in the generating set.  Replacing one of the generators and re-indexing if necessary, we may assume without loss of generality that $\phi_1(\g_1)$ is hyperbolic.  

Then $\phi_2(\g_1)$ is either hyperbolic or parabolic, since it must also be of infinite order (the $\LL\circ\phi_i$'s are isomorphisms onto discrete subgroups of $\soto$).  If it happens to be parabolic, then $\phi_2(\g_1\g_2^k)$ is hyperbolic for large enough $|k|$ and its fixed point set is not invariant under the action of any elliptic generator.  Choose $|k|$ large enough so that $\phi_1(\g_1\g_2^k)$ is also hyperbolic and its fixed point set is not invariant under the action of any elliptic generator.  Thus substituting $\g_1$ for $\g_1\g_2^k$ if necessary, we may assume without loss of generality that both $\phi_1(\g_1)$ and $\phi_2(\g_1)$ are hyperbolic, and their respective fixed point sets are not invariant under the action of any elliptic generator.

Let $j\in\{2,\ldots ,n\}$; then both $\phi_1(\g_1^{k_j}\g_j)$ and  $\phi_2(\g_1^{k_j}\g_j)$ are 
hyperbolic for large enough $|k_j|$.  Substituting if necessary, we may thus assume that 
for $i=1,2$, $\phi_i(\G)$ is freely generated by hyperbolic elements 
$\phi_i(\g_1),\ldots ,\phi_i(\g_n)$.

Next, the following result allows us to further reduce to the case of
a rank two group.  Recall that $G\subset\psl$ is a {\em Schottky}
subgroup of $\psl$ if it admits generators $g_1,\ldots g_n$, called
{\em Schottky 
generators}, such that the following holds: there exist $2n$ disjoint
closed intervals 
$A^\pm _i$ on the boundary of the hyperbolic plane, with 
$g_i(A_i^-)=cl(\partial\hyp\setminus A_i^+)$, $i=1,\ldots,n$, where
$cl$ denotes 
closure.  In particular, if $w=g_{i_1}^{j_1}\cdots g_{i_k}^{j_k}$,
$j_k\in{\mathbb Z}$, is
reduced (i.e. $g_{i_{m+1}}\neq g_{i_m}$, $m=1,\ldots ,k-1$), then $\xp{w}\in A^{\sigma(j_1)}_{i_1}$ and $\xm{w}\in A^{-\sigma(j_k)}_{i_k}$, where $\sigma(j)$ is the sign of $j$.

\begin{prop}\label{thm:Reduction}{\em [Rank-two free subgroups suffice]}~
Let $\G$ be a finitely generated free product of cyclic groups, with
generating set $\g_1,\ldots,\g_n$.  Suppose $\phi_1,\phi_2:\G\longrightarrow\isoaff$ are isomorphisms, such
that each $\LL\circ\phi_i$ is an
isomorphism onto a discrete subgroup of $\psl$ and for $i=1,2$, $\phi_i(\G)$ is freely generated by hyperbolic elements $\phi_i(\g_1),\ldots ,\phi_i(\g_n)$.  
Suppose furthermore that for every rank two subgroup, $H=\langle \g_i,\g_j \rangle\subset \G$, the restriction of
$\phi_1\circ\phi_2^{-1}$ to $H$ is an inner automorphism.   Then
$\phi_1\circ\phi_2^{-1}$ is an inner autmorphism of $\phi_2(\G)$.
\end{prop}

\begin{proof}
The result is clear if the linear parts are the same, so assume that  $\phi_1(\G)$, 
$\phi_2(\G)\subset\soto\cong\psl$.  Also, we assume that $n=3$; the general result is obtained by induction.

Suppose $G=\langle g,h,i\rangle$ and $G'=\langle g',h',i'\rangle$ are discrete subgroups of $\psl$, with hyperbolic 
generators, such that:
\begin{align*}
\langle g',h'\rangle & =\phi\langle g,h\rangle\phi^{-1}\\
\langle g',i'\rangle & =\psi\langle g,i\rangle\psi^{-1}\\
\langle h',i'\rangle & =\pi\langle h,i\rangle\pi^{-1} ,
\end{align*}
for some $\phi,\psi,\pi\in\soto$.   We may further assume that $g,h,i$ (resp. $g',h',i'$) are Schottky 
generators for $G$ (resp. $G'$).  Indeed, $g,h,i$ must have distinct fixed point sets; we get Schottky 
subgroups by substituting for high enough powers of the generators, without affecting $\phi,\psi$ and 
$\pi$.  Then $\psi^{-1}\phi g\phi^{-1}\psi=g$, i.e.,
$\psi^{-1}\phi=\zeta_g$, where $\zeta_g$ commutes with $g$ and thus admits
the same null frame.  Also:
\begin{align*}
\langle g',h'\rangle & =\phi\langle g,h\rangle\phi^{-1} \\
\langle g',i'\rangle & =\phi\zeta_g\langle g,i\rangle\zeta_g^{-1}\phi^{-1} \\
\langle h',i'\rangle & =\phi\zeta_h\langle h,i\rangle\zeta_h^{-1}\phi^{-1} ,
\end{align*}
for some $\zeta_h$ commuting with $h$.  But then $\zeta_g\zeta_h^{-1}$
commutes with $i$, implying that they share the same null frame.  This is
impossible: $\phi_2(\G)$ being a Schottky group, there are pairwise 
disjoint sets
of future-pointing null vectors, one containing the attracting eigenvectors 
of both
$g$ and  $\zeta_g\zeta_h^{-1}$, and the other containing the attracting 
eigenvector
of $i$.  
\end{proof}

Let $\phi=\phi_2\circ\phi_1^{-1}:\phi_1(\G)\rightarrow \phi_2(\G)$.  We 
will use the
following notation: 
\begin{align*}
\Gamma&=\phi_1(\G)=\langle \gg,\gh\rangle\\
\Gamma_\phi&=\phi_2(\G)=\langle
\gg_\phi,\gh_\phi\rangle ,
\end{align*}
where $\gg_\phi=\phi(\gg)$ and $\gh_\phi=\phi(\gh)$.  We will denote the 
linear part and smallest eigenvalue of $\ggp$ (resp. $\ghp$) by $\gp$ and 
$\lagp$ (resp. $\hp$ and $\lahp$).

Since we are interested in conjugacy classes of $\Gamma$, we can
explicitly choose particularly nice representatives of the
conjugacy class in order to simplify our calculations.  For any hyperbolic $g$ and $h$ 
in $\soto$ which generate a non-elementary
group and any three distinct future pointing vectors of unit
Euclidean length $\vv_1,\vv_2,\vv_3$ there is an element
$f\in\oto$  such that $\vx^+\left( fgf^{-1}\right) =\vv_1$,
$\vx^-\left( fgf^{-1}\right) =\vv_2$ and $\vx^-\left(
fhf^{-1}\right) =\vv_3$. This follows directly from the fact that $\oto$ 
acts transitively on triples of distinct points on the boundary of the 
hyperbolic plane, and that:

\begin{equation*}
\frac{f\left( \vx^{\pm} (g)\right) }{ \|
f\left(\vx^{\pm} (g)\right)\| }= \vx^{\pm} (fgf^{-1}).
\end{equation*}

Accordingly, conjugate $\Gamma_\phi$ by $f$ so that:
\begin{align*}
\xp{\gh}&=\xp{\gh_\phi} \\ 
\xm{\gh}&=\xm{\gh_\phi} \\ 
\xm{\gg}&=\xm{\gg_\phi}.
\end{align*}
In particular, the linear part of $h$ is determined up to a choice of 
eigenvalue $\lah$.

Finally, conjugating by a pure
translation translates the invariant line of a Lorentzian
transformation.  We may already assume that the invariant lines $C_\gh$ and 
$C_{\ghp}$ are parallel.  Further conjugate $\Gamma_\phi$ by a pure 
translation taking $C_{\ghp}$ to $C_\gh$.  The translational part of an 
isometry is determined by the displacement factor along its invariant line, 
which is equal to the Margulis invariant.  Thus, if $\ga(\gh)=\ga(\ghp)$, 
the translational parts of $\gh$ and $\ghp$ are equal.

Consequently, it suffices to prove the following restatement of our Main 
Theorem
\ref{thm:Main}.

\begin{thm}\label{thm:Restate}{\em [Technical version of strong isospectrality]}~
Let $\Gamma=\langle \gg, \gh \rangle ,\Gamma_\phi=\langle
\gg_\phi, \gh_\phi \rangle \subset \isoaff$ be 
non-radiant groups with respective linear parts $G_{Id}$ and $G_\phi$, 
such that the generators are all hyperbolic.  Denote by $u$, $u_\phi$, the corresponding cocycles in $Z^1(G_{Id},\rto)$, $Z^1(G_\phi,\rto)$, respectively.  Let $\phi:\Gamma
\rightarrow \Gamma_\phi$  be an isomorphism such that $\phi(\gg)=\gg_\phi$ and
$\phi(\gh)=\gh_\phi$. Suppose that: 
\begin{align*}
\xp{\gh}&=\xp{\gh_\phi} \\ 
\xm{\gh}&=\xm{\gh_\phi} \\ 
C_{\gh}&=C_{\ghp} \\
\xm{\gg}&=\xm{\gg_\phi}
\end{align*}
and furthermore, that the marked Margulis spectra are equal:  
\begin{equation*}
  \ga_{Id}([u])=\ga_{\LL(\phi)}([u_\phi]),
\end{equation*}
where $\Gamma\setminus Id$ and $\Gamma_\phi\setminus Id$ are given compatible orderings. Then $\phi$ is a pure translation.
\end{thm}

The theorem follows from Lemma~\ref{lem:dotsequal} and will thus be proved 
at the end of Section~\ref{sec:main}.

\section{Step two: Estimates for the rate of convergence of the 
eigenvectors}\label{sec:estimates}
We first need to establish a few technical lemmas.

Let $d(\vv, \vw)= \| \vv-\vw \|$ denote the Euclidean distance
between the endpoints of $\vv$ and $\vw$ emanating from the
origin. As $\vx^{\pm}(\gw)$ are normalized to have Euclidean
length $1$, we define $\U$ to be the set of future pointing
Euclidean unit vectors.  For any given $g\in\soto$ we define an
associated map $\tilde{g}:\U\rightarrow\U$ defined as follows:
\begin{equation*} \tilde{g} (\vv) = \frac{g(\vv)}{ \| g(\vv) \| }.
\end{equation*}

\begin{lemma}\label{lemma:close}
Given any hyperbolic element $g\in\soto$ and any $\epsilon >0$,
there exists an $N$ such that for all $n\geq N$ if $\vv\in\U$ such
that $d(\vv, \xm{g})>\epsilon$ then $\tilde{g^n}(\vv)\in\U$
satisfies $d( \tilde{g^n}(\vv) , \xp{g}) < \epsilon$.
\end{lemma}

This lemma is a crude, but particularly clean, version of a
description of the action of a hyperbolic element on null vectors.
For vectors outside the $\epsilon$ neighborhood of $\xm{g}$, we
have the following lemma about how fast the vectors approach
$\xp{g}$.

  \begin{lemma}\label{lemma:closetoxp}
  If $g\in\soto$ is
hyperbolic and $\vv\in\U\setminus \{ \xm{g} \}$, then
$d(\tilde{g}^n(\vv), \xp{g})\sim O(\lambda_{g})$.
\end{lemma}

\begin{proof}
Given any distinct $\xm{g}$, $\xp{g}$, and $\vv$
there is a fixed $h\in\soto$ such that for $g_c = hgh^{-1}$ we
have the following:
  \begin{equation}\label{eq:nicevectors}
\xm{g_c}=\begin{bmatrix} 0\\ \beta \\ \beta\end{bmatrix} ,
\xp{g_c} =\begin{bmatrix} 0\\ -\beta \\ \beta\end{bmatrix}, \vv_c
= \begin{bmatrix} \beta \\ 0\\\beta \end{bmatrix}
\end{equation}
where $\beta= \sqrt{2}/2$ and $\vv_c = \tilde{h}(\vv)$. The map
$\tilde{h}$ changes distances between points on $\U$ by a bounded
multiple. But conjugation by $h$ does not affect the eigenvalues,
i.e. $\gl_{g_c} =\gl_g$.

Therefore, it is enough to show this lemma in the special case
where the vectors are given above. By direct calculation, we have
the following:

\begin{equation}\label{eq:closeorder} d\left( \xp{g_c} ,
\tilde{g_c} ^n(\vv_c) \right) =\lambda^n_{g}
\frac{\sqrt{1+\lambda^{2n}_{g}}}{\beta(1+\lambda_{g}^{2n})}.
\end{equation}
\end{proof}

Now we assert that if the expanding and contracting eigenvectors
of two hyperbolic elements are close, then the fixed vectors are
also close.

\begin{lemma}\label{lemma:closetoxo}
  Given a fixed hyperbolic
element $g\in\soto$  and another hyperbolic element $h\in\soto$,
then

\begin{equation*}
d(\xo{g},\xo{h}) \sim O\left( \max \left(
d(\xm{g},\xm{h}),d(\xp{g},\xp{h}) \right) \right)
\end{equation*}
\end{lemma}
\begin{proof}
It is enough to show this lemma when one pair of expanding or
contracting eigenvectors for $g$ and $h$ are the same. That is,
let
\begin{equation}
\xm{g}=\xm{h}= \begin{bmatrix} 0\\ \beta \\
\beta\end{bmatrix} , \xp{g} =\begin{bmatrix} 0\\ -\beta \\
\beta\end{bmatrix}, \mbox{ and } \xp{h}=\begin{bmatrix}
\beta\sin\delta \\ -\beta\cos\delta \\ \beta\end{bmatrix}.
\end{equation}
By direct calculation
  \begin{equation*}
d(\xp{g},\xp{h}) = \sqrt{1-\cos\delta}\mbox{ and }
d(\xo{g},\xo{h}) =\sqrt{2}|\sin\delta |/ (1+\cos\delta).
\end{equation*}
Then:
\begin{equation*}
\lim_{\delta\rightarrow 0}
\frac{\sqrt{1-\cos\delta} }{(\sqrt{2}|\sin\delta|)/(1+\cos\delta)}
= \lim_{\delta\rightarrow 0} \frac{\beta |\sin\delta|}{\beta|\sin\delta|}=1.
\end{equation*}
\end{proof}

\section{Step three: Equality of the eigenvectors and eigenvalues of the 
generators}\label{sec:main}

At this time, we can attack the main lemma needed to prove
Theorem~\ref{thm:Restate}. Recall that the linear part of $\gg$
(resp. $\gh$, $\ggp$, $\ghp$) is denoted $g$ (resp. $h$, $\gp$,
$\gh$).  Also we are assuming that $\xpm{\ghp}=\xpm{\gh}$ and that
$\xm{\ggp}=\xm{\gg}$.

We need only consider those words of the form $\gh^n\gg^m$ and 
$\ghp^n\ggp^m$, $n,m\geq 0$.  (Similar calculations were done by Goldman~\cite{G}, and the first author~\cite{C}.)  We may assume that these words are all hyperbolic, for large enough $m,n$.  Define:
\begin{equation}\label{eq:xogh}
\xo{g,h}  = \frac{-1}{\ldot{\xm{g},\xp{h}}}\xm{g} \boxtimes\xp{h}.
\end{equation}
This is a unit-spacelike vector to which both $\xo{h^ng^m}$ and 
$\xo{\hp^n\gp^m}$ converge, as $n,m \rightarrow \infty$. In fact, by
Lemma~\ref{lemma:closetoxp} and Lemma~\ref{lemma:closetoxo}, as 
$n,m\rightarrow\infty$, the distance between
$\xo{h^ng^m}$ and $\xo{g,h}$ decreases as $\max \left( \gl_{h}^n,
\gl_{g}^m \right)$, and the distance between $\xo{\hp^n\gp^m}$ and
$\xo{g,h}$ decreases as $\max \left( \gl_{h_\phi}^n,
\gl_{g_\phi}^m \right)$.

\begin{lemma}\label{lem:dotsequal} If
the hypotheses in Theorem~\ref{thm:Restate} are satisfied then
$\xp{\gg}=\xp{\gg_\phi}$.
\end{lemma}

\begin{proof}
For any hyperbolic $\gw\in\isoaff$, we define $E^{\pm}(\gw)$ to be the
plane containing $C_{\gw}$ and parallel to $\xpm{\gw}$.

Let $q = C_{\gg} \cap E^-(\gh)$ and
$p_m = \gg^{-m}(q)$, for $m\geq 0$.  Then $q=r +\gk \xm{h}$, for some $r\in 
C_{\gh}$ and $\gk\in R$.   (See Figure~\ref{fig:points}.)

\begin{figure}
\centerline{\input{Lemma5_1.pstex_t}}
\caption{$q=\cg\cap  E^-(\gh)$ and $p_m = \gg^{-m}(q)$, $m=1,2$.  Then $q=r +\gk \xm{h}$, for some $r\in 
C_{\gh}$ and $\gk\in R$ }
\label{fig:points}
\end{figure}

We can choose similar points for $\ggp$ and $\ghp$; denote by
$\gkp\in\R$ the scalar corresponding to $\gk$.  We will show that
$\gk=\gkp$. To this end, we compute $\ga(\gh^n\gg^m)$:
\begin{equation*}
\ga(\gh^n\gg^m)=\ldot{\gh^n\gg^m(p_m)-p_m,\xo{h^ng^m}}.
\end{equation*}

Write $\gh^n\gg^m(p_m)-p_m$ as $\gh^n(q)-q+q-p_m$.  Then:
\begin{align*}
q-p_m &=m\ga(\gg)\xo{g} \\
\gh^n(q)-q &=n\ga(\gh)\xo{h} +\gk(\lah^n-1)\xm{h}.
\end{align*}
Since $\xo{h^ng^m}$ converges to $\xo{g,h}$ faster than $n,m$, 
$\ga(\gh^n\gg^m)$ asymptotically approaches:
\begin{equation*}
m\ga(\gg)\ldot{\xo{g},\xo{g,h}}+n\ga(\gh)\ldot{\xo{h},\xo{g,h}}+\gk(\lah^n-1)\ldot{\xm{h},\xo{g,h}}.
\end{equation*}
By Lemma~\ref{lem:xpxo} and Equation~\eqref{eq:xogh}:
\begin{align*}
\ldot{\xo{g},\xo{g,h}}&= \frac{-1}{\ldot{\xm{g},\xp{h}}}\ldot{\xo{g},\xm{g} 
\boxtimes\xp{h}} \\
                                        &=\frac{-1}{\ldot{\xm{g},\xp{h}}}\ldot{\xp{h},\xo{g} 
\boxtimes\xm{g}} \\
                                        &=\frac{1}{\ldot{\xm{g},\xp{h}}}\ldot{\xp{h},\xm{g}}=1.
\end{align*}
Thus $\ga(\gh^n\gg^m)$ asymptotically approaches:
\begin{equation*}
m\ga(\gg)+n\ga(\gh)+\gk(\lah^n-1)\ldot{\xm{h},\xo{g,h}},
\end{equation*}
and in the same fashion, $\ga(\ghp^n\ggp^m)$ asymptotically approaches:
\begin{equation*}
m\ga(\ggp)+n\ga(\ghp)+\gkp(\lambda_{\eta_\phi}^n-1)\ldot{\xm{\hp},\xo{g,h}}.
\end{equation*}

Since $\ga(\ghp^n\ggp^m)=\ga(\gh^n\gg^m)$ for all $n,m\geq 0$, and since 
$\xm{\hp}=\xm{h}$, it follows that $\gk=\gkp$, as claimed.

We have assumed that $C_{\ghp}=C_\gh$; thus $C_{\gg}$ and $C_{\gg_\phi}$ 
must intersect the
line  $C_{\gh} -\gk \xm{h}$.

By similar reasoning, replacing $\gh$
and $\gh_\phi$ with their inverses above, we know that  $C_{\gg}$
and $C_{\gg_\phi}$ must also intersect the line  $C_{\gh} -
\gk_{-}\xp{h}$, where $\gk_{-}\xp{h}$ is a vector which points
from a point on $C_{\gg}$ to a point on $C_{\gh}$. That is,
$C_{\gg}$ and $C_{\gg_\phi}$ must both lie in the plane defined by
the disjoint parallel lines $C_{\gh} -\gk \xm{h}$ and $C_{\gh}
-\gk_{-} \xp{h}$.

We have assumed that $\xm{g}=\xm{g_\phi}$. The planes $E_{\gg}^-$ and
$E_{\gg_\phi}^-$ are Lorentzian perpendicular to  $\xm{g}$ and
$\xm{g_\phi}$, so $E_{\gg}^-$ and $E_{\gg_\phi}^-$ are parallel.
The lines $C_{\gg}$ and $C_{\gg_\phi}$ are the intersections of
two parallel planes with one fixed plane, so they must be
parallel. Alternatively, $\xo{\gg}=\xo{\gg_\phi}$ so
$\xp{\gg}=\xp{\gg_\phi}$.

\end{proof}

Lemma 5.1 can be restated as follows:

\begin{lemma}{\em [Fixed point isospectrality]}~
Consider a discrete $G\subset 
Isom(\hyp)$  which is  purely hyperbolic, i.e. all non-identity elements are 
hyperbolic. Let:
\begin{eqnarray*}
f: G & \longrightarrow \partial \hyp\times \partial \hyp \\
   g & \longmapsto (\xp{g},\xm{g}),
\end{eqnarray*}
where $\xp{g}$ (resp. $\xm{g}$) is the attracting (resp. repelling) fixed point of $g$.  Then $f$ completely determines the group $G$.
\end{lemma}

\begin{proof} (Theorem~\ref{thm:Restate})

The hypotheses and Lemma~\ref{lem:dotsequal} imply that the
eigenspaces of $g$ and $\gp$ are the same, as are those of $h$ and
$\hp$.    Also, by hypothesis, $\vv_h=\vv_{\hp}$.

The equalities $\gl_{\gg}=\gl_{\gg_\phi}$ and
$\gl_{\gh}=\gl_{\gh_\phi}$  follow immediately from
Lemma~\ref{lem:dotsequal}.  Indeed, substituting $ghg^{-1}$ for
$h$, we obtain $\xp{ghg^{-1}}=\xp{g_\phi h_\phi g_\phi^{-1}}$.
The fact that $\xp{ghg^{-1}}$ is parallel to $g(\xp{h})$ uniquely
determines $\gl_{\gg}$.  We show that $\lah=\lahp$ in a similar
fashion, by considering $hgh^{-1}$.

Finally, we consider the translational parts of $\gg$ and $\ggp$.
The proof of Lemma~\ref{lem:dotsequal} shows that $C_\gg$ and
$C_{\ggp}$ are parallel; as a matter of fact, since both lines
intersect $C_\gh+\gk\xm{h}$, $C_{\ggp}$ can be translated to
$C_\gg$ by a translation $\tau$, parallel to $\xo{h}$.  This
conjugation preserves $C_h$, thus $\gh$.  Therefore:
\begin{equation*}
\Gamma = \tau\Gamma_\phi\tau^{-1},
\end{equation*}
where $\tau$ is a pure translation.
\end{proof}


\end{document}

%% file: Lemma5_1.pstex_t
\begin{picture}(0,0)%
\epsfig{file=Lemma5_1.pstex}%
\end{picture}%
\setlength{\unitlength}{0.00062500in}%
\begingroup\makeatletter\ifx\SetFigFont\undefined%
\gdef\SetFigFont#1#2#3#4#5{%
  \reset@font\fontsize{#1}{#2pt}%
  \fontfamily{#3}\fontseries{#4}\fontshape{#5}%
  \selectfont}%
\fi\endgroup%
\begin{picture}(5064,5481)(1174,-5830)
\put(4516,-2536){\makebox(0,0)[lb]{\smash{\SetFigFont{9}{10.8}{\rmdefault}{\mddefault}{\updefault}$r$}}}
\put(3526,-3511){\makebox(0,0)[lb]{\smash{\SetFigFont{9}{10.8}{\rmdefault}{\mddefault}{\updefault}$\kappa\xm{h}$}}}
\put(5821,-4006){\makebox(0,0)[lb]{\smash{\SetFigFont{9}{10.8}{\rmdefault}{\mddefault}{\updefault}$p_2$}}}
\put(5116,-4336){\makebox(0,0)[lb]{\smash{\SetFigFont{9}{10.8}{\rmdefault}{\mddefault}{\updefault}$p_1$}}}
\put(4261,-4681){\makebox(0,0)[lb]{\smash{\SetFigFont{9}{10.8}{\rmdefault}{\mddefault}{\updefault}$q$}}}
\put(1411,-5791){\makebox(0,0)[lb]{\smash{\SetFigFont{9}{10.8}{\rmdefault}{\mddefault}{\updefault}$\cg$}}}
\put(6106,-2851){\makebox(0,0)[lb]{\smash{\SetFigFont{9}{10.8}{\rmdefault}{\mddefault}{\updefault}$C_\eta$}}}
\end{picture}